\documentclass[letterpaper, 10 pt, conference]{ieeeconf}  

\IEEEoverridecommandlockouts                              
\usepackage{graphics} 
\usepackage{epsfig} 
\usepackage{times} 

\usepackage{ntheorem}
\usepackage{amsmath,amssymb,amsfonts}
\usepackage{xcolor}

\usepackage[compress]{cite}

\newtheorem{lemma}{Lemma}[section]
\newtheorem{definition}{Definition}[section]
\newtheorem{theorem}{Theorem}
\newtheorem{assumption}{Assumption}

\title{\LARGE \bf
Data-Driven Performance Guarantees for\\ Parametric Optimization Problems}

\author{Jingyi Huang, Paul Goulart and Kostas Margellos
\thanks{The authors are with the Department of Engineering Science,
        University of Oxford, Oxford, United Kingdom. E-mails:
        \{jingyi.huang, 
        paul.goulart,  kostas.margellos\}@eng.ox.ac.uk}%
\thanks{For the purpose of Open Access, the authors have applied a CC BY copyright license to any Author Accepted Manuscript (AAM) version arising from this submission.}
}

\begin{document}

\maketitle
\thispagestyle{empty}
\pagestyle{empty}

\begin{abstract}
We propose a data-driven method to establish probabilistic performance guarantees for parametric optimization problems solved via iterative algorithms. Our approach addresses two key challenges: providing convergence guarantees to characterize the worst-case number of iterations required to achieve a predefined tolerance, and upper bounding a performance metric after a fixed number of iterations. These guarantees are particularly useful for online optimization problems with limited computational time, where existing performance guarantees are often unavailable or unduly conservative. We formulate the convergence analysis problem as a scenario optimization program based on a finite set of sampled parameter instances. Leveraging tools from scenario optimization theory enables us to derive probabilistic guarantees on the number of iterations needed to meet a given tolerance level. Using recent advancements in scenario optimization, we further introduce a relaxation approach to trade the number of iterations against the risk of violating convergence criteria thresholds. Additionally, we analyze the trade-off between solution accuracy and time efficiency for fixed-iteration optimization problems by casting them into scenario optimization programs. Numerical simulations demonstrate the efficacy of our approach in providing reliable probabilistic convergence guarantees and evaluating the trade-off between solution accuracy and computational cost.
\end{abstract}

\section{INTRODUCTION}
Online parametric optimization problems  arise frequently in control theory, with model predictive control (MPC) \cite{MPC_book,maciejowski2007predictive, Jingyi_MPC} being a key example in which the parameters correspond to the initial states. A major challenge in online optimization is the limited computational time available, often restricting the solver to a fixed number of iterations. In online optimization with strict time constraints such as MPC, analytical convergence guarantees are often very conservative relative to the actual optimization performance observed in practice. This raises two critical questions: (i) Can we ensure convergence within a given iteration budget based on numerical experiment?; and (ii) How can we determine the minimal number of iterations needed to achieve satisfactory performance? In this paper, we address these questions by providing a probabilistic convergence guarantee, learning the minimum number of iterations required for convergence, and quantifying the trade-off between iterations and performance.

Throughout, we will consider parametric constrained optimization problems of the form, 
\begin{equation}\label{parametric_optimization_problem}
\begin{split}
    \min_{z} \quad & f(z, x) \\
    \text{subject to:} \quad & f_i(z, x) = 0, \quad i = 1, \ldots, I,\\
    & g_j(z, x) \geq 0, \quad j = 1, \ldots, J,
\end{split}
\end{equation}
where \( z \in \mathbb{R}^n \) is the optimization variable, \( x \in \Delta \subset \mathbb{R}^d \) is the problem parameter, and \( \Delta \) is the parameter space representing the set of all possible values of \( x \) for which the parametric optimization problem has at least one feasible solution for \( z \). The (potentially local) optimal solution to the problem \eqref{parametric_optimization_problem} will be denoted $z^*(x)$.  We will assume that such a solution is unique (for any given parameter); otherwise a particular solution can be singled-out by means of a tie-break rule.
Problems in the form of \eqref{parametric_optimization_problem} typically appear in optimization-based control techniques such as MPC, where they must be solved repeatedly  for varying initial conditions \( x \). In this context, the optimization variable \( z \) represents the open-loop state and input sequences, while the objective function \( f: \mathbb{R}^n \times \mathbb{R}^d \rightarrow \mathbb{R} \) defines the criterion to be minimized for optimal control. The functions $f_i(z, x)$ and $g_j(z, x)$ account for the system dynamics, and the input and state constraints, respectively.

Theoretical convergence analysis for iterative optimization methods generally provide a-priori convergence rates based only on the problem structure. 
An overview of convergence analysis across different optimization methods can be found in \cite[Chapters 9-15]{Beck2017} and \cite{RyuYin2022}. Some studies have investigated worst-case convergence rates for unconstrained optimization problems via computer-assisted analysis \cite{DasGupta2024, Taylor2017}. However, since these approaches do not exploit the properties of the problem parameters, their convergence rates tend to be conservative. Incorporating the characteristics of the problem parameters could lead to a tighter convergence guarantee. 

The authors of \cite{sambharya_data-driven_2024} proposed a data-driven approach based on sample convergence bounds \cite{Langford2001NotBT, JMLR:v6:langford05a}. Their method provides probabilistic performance guarantees for user-specified iteration counts, ensuring each evaluation holds with a pre-defined confidence level of \(1 - \beta\), where $\beta \in (0,1)$. By evaluating \(N^{\mathrm{eval}}\) distinct iteration counts, they quantify the trade-off between optimization performance and iteration count, with an overall confidence level of \(1 - N^{\mathrm{eval}} \beta\). However, this method cannot leverage sample-dependent information to select iteration counts for performance evaluation, potentially requiring more evaluations to establish an informative trade-off between computation time and performance, thereby reducing the overall confidence level.

Inspired by \cite{sambharya_data-driven_2024}, we propose a data-driven method based on the so called scenario approach \cite{modulating_robustness_2013, GarattiCampi2022, Marco_CDC_2020, wait_and_judge, kostas_2014, Kostas_2015, non-convex_scenario_optimization2024} to establish sharp convergence guarantees with high confidence. Our approach solves a scenario optimization problem to automatically determine the worst-case iteration count required to achieve a prescribed tolerance level across sampled parameters, while providing probabilistic guarantees that certify its validity for unseen instances. As such, our approach integrates sample information through scenario programming, enabling us to simultaneously \emph{learn} the iteration count for convergence and provide a \emph{generalization guarantee}.

Our contributions can be summarized as follows:
\begin{itemize}
\item We propose a data-driven method to establish probabilistic convergence guarantees for parametric optimization problems.
\item We construct data-driven guarantees for the upper bound of a user-defined performance metric in parametric optimization problems.
\item We demonstrate the effectiveness of our method using an MPC example, showcasing its ability to provide probabilistic convergence guarantees and evaluate the trade-off between iteration count and performance.
\end{itemize}

Section II formulates the convergence guarantee analysis as robust and relaxed scenario optimization problems, providing probabilistic guarantees for each formulation. Section III provides probabilistic guarantees for the upper bound of a user-defined performance metric in parametric optimization problems. Section IV illustrates the effectiveness of our method through numerical examples. Finally, Section V concludes the paper.


\section{Data-driven convergence guarantees}\label{section:convergence guarantee}
We consider optimization problems of the form \eqref{parametric_optimization_problem}, focusing on the convergence of iterates $z^k(x)$ to an optimizer $z^*(x)$, i.e., $\lim_{k \to \infty} \| z^k(x) - z^*(x) \| = 0$, where $k$ denotes the iteration count during the optimization process. Both the (local) optimal solution \( z^*(x) \) and the iterative optimization variables \( z^k(x) \) depend on the problem parameters \( x \). We denote the set of local optimizers for problem \eqref{parametric_optimization_problem} by \(\mathcal{Z}_{\text{local}}^*(x)\). By our definition of \(\Delta\), the set \(\mathcal{Z}_{\text{local}}^*(x)\) is non-empty. In practice, \( z^*(x) \in \mathcal{Z}_{\text{local}}^*(x) \) would be the solution returned by a solver once a given numerical tolerance is met.
%

In this section, we formulate the problem of establishing convergence guarantees as a scenario optimization problem, learning the number of iterations for convergence and providing a generalization guarantee through the scenario approach.

\subsection{Problem formulation}

We define \( n(x) \) as the number of iterations required for an optimization algorithm to solve problem \eqref{parametric_optimization_problem} with parameter \( x \in \Delta\) to within a specified tolerance. In other words, for a given $x$, $\|z^{n(x)}(x) - z^*(x)\| \leq \eta$, where $\eta$ is a prespecified, user-chosen tolerance level.
We denote by \( n^* \) the minimum number of iterations needed so that this tolerance is met for all possible realizations of $x$. 
Determining \( n^* \) provides an upper bound on the computational effort required to solve the parametric optimization problem within the desired accuracy. It can be mathematically formulated as

\begin{equation}\label{n robust}
    n^* = \sup_{x \in \Delta} \;\; n(x).
\end{equation}
%
%

The parameter space \( \Delta \) comprises all parameter values \( x \) for which the optimization problem \eqref{parametric_optimization_problem} has at least one feasible solution. Even if \( \Delta \) is known, the exact value of \( n^* \) cannot be computed in general as \( \Delta \) could be a continuous set. To circumvent this, we employ a data-driven approach to estimate \( n^* \). Specifically, we assume that \( \Delta \) is endowed with a $\sigma$-algebra and associated probability measure $\mathbb{P}$. We then consider the following \emph{scenario program}: 
\begin{equation}\label{n scenario}
\begin{split}
    n^*_N = \, y^*_N =   \min_{y\geq0} \; &y\\
    \text{subject to:} \quad  &y \geq \,n(x_i), \quad \forall i = 1, \ldots, N,
\end{split}
\end{equation}
where \( N \) denotes the number of sampled parameter instances \(\{x_i\}_{i=1}^N\), with each \( x_i \) drawn independently according to \( \mathbb{P} \). Thus, \( n^*_N \) provides an empirical estimate for \( n^* \).

As \( n^*_N \) is an empirical estimate based on observed data, it is essential to establish a generalization guarantee to ensure its validity on unseen instances. This is achieved via the scenario approach, which provides probabilistic generalization guarantees.


\subsection{Probabilistic convergence guarantees}
In problem \eqref{n scenario}, each sample $x_i$ defines a constraint of the form $y \geq n(x_i)$, which can be written as $y \in \mathcal{Y}_{x_i}$, where $\mathcal{Y}_{x_i}$ denotes the set of all feasible values of $y$ for the $i$-th scenario. To evaluate how well the solution $n^*_N$ generalizes to new samples $x$, we quantify the \emph{risk} associated with $n^*_N$.

\begin{definition}[Risk] \label{def:risk}
The risk of a solution \( y \) is the probability that it violates the scenario constraint induced by a new sample:
\begin{equation}
V(y) = \mathbb{P}\{x \in \Delta :~ y \notin \mathcal{Y}_{x}\}.
\end{equation}
\end{definition}
In our context, \(y^*_N = n_N^*\), and \( V(y^*_N) = \mathbb{P}\{x \in \Delta :~ n(x) > n_N^*\}\) denotes the probability that a new parameter instance \( x \) requires more iterations to solve than the robust estimator \( n_N^* \).

\begin{theorem}[{\protect{\cite[Theorem 1]{campi2008exact}}}]\label{theorem: robust} Consider a convex scenario optimization problem with optimization variable $y \in \mathbb{R}^d$. Let $\{x_i\}_{i=1}^N$ be $N$ i.i.d. samples from $(\Delta, \mathcal{F}, \mathbb{P})$, where $\mathcal{F}$ is the event space, and $y_N^*$ be the optimal solution (assumed to exist and to be unique). Fix $\beta \in (0,1)$. We have that
\begin{equation}\label{eq:theorem 1}
\mathbb{P}^N\{\{x_i\}_{i=1}^N:~V(y_N^*) \leq \epsilon\} \geq 1 - \beta,
\end{equation}
where $\epsilon \in (0,1)$ is such that $\beta := \sum_{i=0}^{d-1} \binom{N}{i} \epsilon^{i} (1-\epsilon)^{N-i}$.
\end{theorem}

With confidence at least $1-\beta$, the quantity $\epsilon$ is an upper bound on the risk $V(y_N^*)$. 
Since problem \eqref{n scenario} is convex and admits a unique minimizer, Theorem \ref{theorem: robust} can provide a generalization guarantee for the empirical estimate $n^*_N$. Given that the problem in \eqref{n scenario} is a convex scenario program and involves a scalar decision variable $d = 1$, Theorem \ref{theorem: robust} is directly applicable and results in $\epsilon = 1 - \beta^{\frac{1}{N}}$. Therefore, we obtain the probabilistic guarantee that, with confidence at least $1-\beta$, for any new parameter realization $x$, the probability of requiring more than $n^*_N$ iterations to solve is bounded above by $\epsilon$.


\subsection{Probabilistic convergence guarantees via constraint relaxation}

There is an inherent trade-off between the allocated computation time and the probability of successfully solving the optimization problem within this time budget. Previously, we framed the problem of offering probabilistic convergence guarantees as a scenario optimization problem, yielding a solution such that for all parameter realizations/samples the number of iterations needed to meet a given tolerance level is at most equal to $n_N^*$.
However, real-time applications may impose limits on computational time. 

To trade off available computation time against the risk of meeting a given convergence tolerance, we relax the constraints in \eqref{n scenario}, leading to the following problem: 
\begin{equation}\label{n_relaxed scenario optimization}
\begin{split}
    \min_{y \geq 0, \, \xi_{N,i} \geq 0, \, i=1,\ldots,N} \quad & y + \rho \sum_{i=1}^N \xi_{N,i}\\
    \text{subject to:} \quad & y + \xi_{N,i} \geq n(x_i), \quad i = 1, \ldots, N,
\end{split}
\end{equation}  
where \( \xi_{N,i}, i = 1, \ldots, N \) are relaxation variables. When \( \xi_{N,i} > 0 \), the original constraint \( y \geq n(x_i) \) is relaxed to \( y + \xi_{N,i} \geq n(x_i) \). In analogy with the results of the previous section, we denote the optimal solution \( y_N^* \)\footnote{We will assume that such an optimal solution \( y_N^* \) is unique; otherwise a particular solution can be singled-out by means of a tie-break rule.}
as \( n^*_N(\rho) \), which now depends on the weight \( \rho > 0 \), which balances the trade-off between computation time and the risk associated with meeting a given convergence tolerance. We denote the number of ``violated'' constraints for the original problem by $q_N^* = |\{i : \xi_{N,i}^* > 0\}|$, which serves as a key metric for evaluating the solution's robustness. Notice, that $q_N^*$ is an \emph{a posteriori} quantity as it depends  on the observed samples. As such, the risk assessments provided in the sequel will be \emph{a posteriori} as well.

Most often in the scenario approach literature with \emph{a posteriori} assessments, a \emph{non-accumulation} assumption on the system constraints \cite{GarattiCampi2022} is imposed. If applied to problem \eqref{n_relaxed scenario optimization}, this would require that for every \( y \geq 0 \), \(\mathbb{P}\{x \in \Delta: y = n(x)\} = 0\). This assumption implies that the scenario constraints \( y \geq n(x_i) \) do not accumulate with non-trivial probability mass at the constraint boundary \( y = n(x) \). However, this assumption naturally does not hold in our setting. While scenarios are drawn from a continuous space, the function \( n(\cdot) \) maps each scenario to a discrete, finite set of values. Thus, coincident constraints, i.e., \( n(x_i) = n(x_j) \) for some \( i \neq j \), occur with non-zero probability. Consequently, classical convex scenario approaches \cite{GarattiCampi2022, Marco_CDC_2020, modulating_robustness_2013} are not applicable to our problem. Recent advances in the scenario approach \cite{non-convex_scenario_optimization2024} allow the characterization of the risk level in problems with constraint relaxation, even if constraints accumulate. This enables us to obtain results similar to that of Theorem \ref{theorem: robust}, however, under some different technical conditions which we introduce below. 

Problem \eqref{n_relaxed scenario optimization} is a decision-making process in which we introduce an augmented decision variable \( \tilde{y} = (y, q) \in \mathbb{R}^2 \). This variable $\tilde{y}$ simultaneously captures \( y \) and the number of relaxed constraints \( q\). The set of all feasible values of $\tilde{y}$ for the $i$-th scenario is denoted by $\tilde{\mathcal{Y}}_{x_i}$. The process can be formalized as a mapping \( M_N^{\rm relax} : \Delta^N \rightarrow \mathbb{R}^2 \), which takes \( N \) i.i.d. samples \( \{x_i\}_{i=1}^N \) from the sample space \( \Delta \) and returns the decision \( \tilde{y}_N^* = (y_N^*, q_N^*) \). To emphasize the dependence on the samples, we also use the notation \( \tilde{y}_N^* = M_N^{\rm relax}(x_1, \ldots, x_N) \).  

To establish generalization guarantees via \cite{non-convex_scenario_optimization2024}, we must verify that \( M_N^{\rm relax} \) satisfies the \emph{consistency} property.

\begin{definition}[Consistency property, {\protect{\cite[Property 1]{non-convex_scenario_optimization2024}}}]\label{consistency_property}
Consider a map $M_m: \Delta^m \rightarrow \tilde{\mathcal{Y}}$ that generates a decision vector $\tilde{y}^*_m$. For any integers \( N \geq 0 \) and \( K > 0 \) and for any choice of \(x_1, \ldots, x_N\) and \(x_{N+1}, \ldots, x_{N+K}\), the following three conditions hold:
\begin{enumerate}
    \item[(i)] If \(x_{i_1}, \ldots, x_{i_N}\) is a permutation of \(x_1, \ldots, x_N\), then it holds that $ M_N(x_1, \ldots, x_N) = M_N(x_{i_1}, \ldots, x_{i_N});$
    \item[(ii)] If \(\tilde{y}^*_N \in \tilde{\mathcal{Y}}_{x_{N+i}}\) for all \(i = 1, \ldots, K\), then it holds that $\tilde{y}^*_{N+K} = M_{N+K}(x_1, \ldots, x_{N+K}) = M_N(x_1, \ldots, x_N) = \tilde{y}^*_N;$
    \item[(iii)] If \(\tilde{y}^*_N \notin \tilde{\mathcal{Y}}_{x_{N+i}}\) for one or more \(i \in \{1, \ldots, K\}\), then it holds that $\tilde{y}^*_{N+K} = M_{N+K}(x_1, \ldots, x_{N+K}) \neq M_N(x_1, \ldots, x_N) = \tilde{y}^*_N.$
\end{enumerate}
\end{definition}

Condition (i) means that the decision of the map \( M_m \) remains unchanged regardless of the order of scenarios. Condition (ii) ensures that the  decision remains the same when additional scenarios are introduced, provided that the previous decision already satisfies the constraints imposed by these new scenarios. Condition (iii) requires that the decision-making process $M_m$ reacts when additional scenarios are added for which the previous decision does not meet the constraints induced by these new scenarios.

\begin{lemma}\label{relaxed_consistency_proof}
    Consider the relaxed scenario optimization problem \eqref{n_relaxed scenario optimization}. The mapping \( M_N^{\rm relax}: \Delta^N \rightarrow \mathbb{R}^2 \) satisfies the consistency property in Definition \ref{consistency_property}.
\end{lemma}

\textit{Proof.} Condition (i) holds because $y_N^*$ and $q_N^*$ are independent of the order of constraints, so is $\tilde{y}_N^*$. Now, examine conditions (ii) and (iii). Add $K$ new scenarios $x_{N+1}, \ldots, x_{N+K}$ to the original samples $x_1, \ldots, x_N$ and assume that \( y_N^* \geq n(x_{N+i}) \) for all \( i = 1, \ldots, K \). Consider problem \eqref{n_relaxed scenario optimization} with $N+K$ in place of $N$. Since \( y_N^* \geq n(x_{N+i}) \) for all \( i = 1, \ldots, K \) meaning $ (y_N^*, \xi_{N,1}^*, \ldots, \xi_{N,N}^*, 0, \ldots, 0)$ is feasible for problem \eqref{n_relaxed scenario optimization} with $N+K$ scenarios. We claim that this is indeed the optimal solution, which can be shown by contradiction. If the optimal solution were different, and it is denoted \( (\bar{y}, \bar{\xi}_{N+ K,i}) \) for \( i = 1, \dots, N+K \), then one of the following two cases must hold:

\begin{enumerate}
    \item[(a)] \( \bar{y} + \rho \sum_{i=1}^{N+K} \bar{\xi}_{N+ K,i} < y^*_N + \rho \sum_{i=1}^{N} \xi_{N,i}^* \). This implies that \( \bar{y} + \rho \sum_{i=1}^{N} \bar{\xi}_{N+ K,i} < y^*_N + \rho \sum_{i=1}^{N} \xi_{N,i}^*\) since the omitted terms \(\bar{\xi}_{N+ K,i} \geq 0, \text{ for } i = N+1, \ldots, N+K\). This then implies that in problem \eqref{n_relaxed scenario optimization}, the solution \( \{\bar{y}, \bar{\xi}_{N+ K,i}\}, i = 1, \dots, N \) achieves a better objective value than the optimal solution \( \{y_N^*, \xi^*_{N,i}\}, i = 1, \dots, N \), leading to a contradiction.

    \item[(b)] \( \bar{y} + \rho \sum_{i=1}^{N+K} \bar{\xi}_{N+ K,i} = y^*_N + \rho \sum_{i=1}^{N} \xi_{N,i}^*\) and $\bar{y} \neq y^*_N$. By omitting terms \(\bar{\xi}_{N+ K,i} \geq 0, \text{ for } i = (N+1), \ldots, (N+K)\), we have \(\bar{y} + \rho \sum_{i=1}^{N} \bar{\xi}_{N+ K,i} \leq y^*_N + \rho \sum_{i=1}^{N} \xi_{N,i}^*.\) If this inequality holds strictly, it would lead to a contradiction similar to case (a). On the other hand, if equality holds, it implies that problem \eqref{n_relaxed scenario optimization} with $N$ scenarios has multiple optimal solutions and \( (\bar{y}, \bar{\xi}_{N+ K,i}) \) for \( i = 1, \dots, N \) ranks better than $ (y_N^*, \xi_{N,1}^*, \ldots, \xi_{N,N}^*)$. This contradicts the tie-break rule, which uniquely selects $ (y_N^*, \xi_{N,1}^*, \ldots, \xi_{N,N}^*)$ as the optimal solution for problem \eqref{n_relaxed scenario optimization} with $N$ scenarios. 
\end{enumerate}

Hence, we have proven that $(y_N^*, \xi_{N,1}^*, \ldots, \xi_{N,N}^*, 0, \ldots, 0)$ is the unique optimal solution for problem \eqref{n_relaxed scenario optimization} when extended to \( N+K \) constraints. Consequently, the decision returned by the mapping \( M_N^{\text{relax}} \) remains unchanged, i.e., $\tilde{y}_{N+K}^* = (y_{N+K}^*, q^*_{N+K}) = (y_{N}^*, q^*_{N}) = \tilde{y}_{N}^*$, which proves condition (ii) in the consistency property.

Now, suppose instead that \( y_N^* < n(x_{N+i}) \) for some \( i \in \{1, \ldots, K\} \). If \( y_{N+K}^* = y_N^* \), then $ \xi_{N+ K, i}^* = \xi_{N, i}^* \text{ for } i = 1, \ldots, N \text{ and } \xi_{N+ K, N+i}^* > 0 \text{ for some } i \in \{1, \ldots, K\}$. As a result, $q_{N+K}^* > q_N^*,$ which implies that \(\tilde{y}_{N+K}^* \neq \tilde{y}_{N}^*\). Alternatively, if \(y_{N+K}^* \neq y_N^*\), it directly follows that \( \tilde{y}_{N+K}^* \neq \tilde{y}_N^* \). This proves condition (iii). \hfill \(\square\)

Having established the consistency property of the mapping \( M_N^{\rm relax} \), we now leverage the results from \cite{non-convex_scenario_optimization2024} to derive generalization guarantees for problem \eqref{n_relaxed scenario optimization}. We first present the key equations that relate the confidence level \(\beta\) to the violation probability \(\epsilon\) \cite{non-convex_scenario_optimization2024}. For a user-specified confidence level \(\beta \in (0,1)\), consider the following polynomial equation in \( t \) for any \( q = 0, 1, \ldots, N-1 \):
\begin{equation}\label{compute_t}
\frac{\beta}{N} \sum_{m=q}^{N-1} \binom{m}{q} t^{m-q} - \binom{N}{q} t^{N-q} = 0,
\end{equation}
where \( N \) represents the number of samples. For each \( q \),~\eqref{compute_t} admits a unique solution \( t(q) \) in the interval \( (0, 1) \). The corresponding violation probability \( \epsilon(q) \) is then given by:
\begin{equation}\label{compute_epsilon}
\epsilon(q) := 
\begin{cases} 
1 - t(q), & \text{for } q = 0, 1, \ldots, N-1, \\
1, & \text{for } q = N.
\end{cases}
\end{equation}

\begin{theorem} Consider the relaxed scenario optimization problem \eqref{n_relaxed scenario optimization}. Fix $\beta \in (0,1)$. We then have
\begin{equation}\label{relaxed_pac}
\mathbb{P}^N \{ \{x_i\}_{i=1}^N:~ V(n^*_N(\rho)) \leq \epsilon(s^*_N)\} \geq 1 - \beta,
\end{equation}
where $\epsilon(q), q = 0,\ldots,N$ is defined as in \eqref{compute_epsilon}, and $s_N^* = 1 + q_N^*$.
\end{theorem}

\textit{Proof.} Lemma \ref{relaxed_consistency_proof} shows that the mapping \( M_N^{\rm relax} \) satisfies the consistency property, which permits the use of \cite[Theorem 10]{non-convex_scenario_optimization2024} to bound the risk $V(y^*_N)$ for solutions of problem \eqref{n_relaxed scenario optimization} at confidence level \(\beta\). The parameter \( s^*_N \) characterizes the solution complexity of \( y^*_N \) through two components: (i) the number of scenarios \( x_i \) with nonzero relaxations (\(\xi^*_{N,i} > 0\)), denoted by \( q^*_N \); and (ii) the cardinality of a minimal set of additional $x_i$ required to uniquely determine \( y^*_N \), which is $1$, since the original problem has only one decision variable and is convex. Therefore, we use $s^*_N = 1 + q^*_N$. \hfill \(\square\)

The regularization parameter \(\rho\) governs the trade-off between minimizing the computation time and penalizing for the violation cost. We next suggest an approach to selecting a suitable value of \(\rho\). Problem \eqref{n_relaxed scenario optimization} can be formulated as an optimization problem $\min_{y \geq 0}\; f(y)$, where

\begin{equation}
    f(y) =  y + \rho \sum_{i=1}^N \max\{0, n(x_i) - y\}. 
\end{equation}

The objective function has a subgradient $1 - \rho \sum^N_{i=1} h_i(y) \in \partial f(y)$, where  

\begin{equation}\label{define_h}
\begin{split}
    h_i(y) := \left\{
\begin{array}{ll}
      0, & \text{if } y \geq n(x_i), \\
      1, & \text{if } y < n(x_i).
\end{array}
\right. \\
\end{split}
\end{equation}

We focus on cases where the constraint \( y \geq 0 \) remains inactive at optimality, as the solution \( y^*_N = 0 \) gives a trivial violation probability \( V(y^*_N) = 1 \). Thus, at optimality, the stationarity condition yields: 
\[0 = 1 - \rho \sum^N_{i=1} h_i(n^*_N(\rho)) \in \partial f(y)|_{y = y^*_N = n^*_N(\rho)}.\] 
Solving for \(\rho\) gives:

\[
\rho = \frac{1}{\sum^N_{i=1} h_i(n^*_N(\rho))}.
\]

The denominator counts the number of violated constraints \(y^*_N < n(x_i)\), i.e., $q^*_N$. While this quantity is only known post-optimization, we can approximate it using a predefined target violation number \(\hat{q}_N\), which is sample-independent. This computation allows for more informed practical choices of the parameter vector, e.g., 
\[
\rho = \frac{1}{\hat{q}_N},
\]
where \(\hat{q}_N\) serves as a proxy for the actual number of violated constraints \(q^*_N\) \footnote{
Following the primal-dual analysis in \cite{New_SVM}, $\hat{q}_N$ can serve as an upper bound for $q^*_N$. This allows the violated constraints to be counted a priori, enabling the use of Theorem 1. The corresponding results are plotted in Figure \ref{fig:n trade off} as red crosses.}.

\textbf{Remark 1.} In this section we have developed a scenario approach framework that establishes generalization guarantees for both \( V(n^*_N) \) and \( V(n^*_N(\rho)) \), where \( n^*_N \) and \( n^*_N(\rho) \) are learned directly from data through scenario optimization. In contrast, \cite{sambharya_data-driven_2024} proposes a method to bound the risk \( V(n_a) \) for predefined, sample-independent iteration numbers \( n_a \). Their methodology suffers from two critical limitations: first, obtaining probabilistic convergence guarantees requires evaluating numerous candidate \( n_a \) values, with confidence degrading linearly in the number of evaluations; second, the inability to incorporate sample information during iteration number selection leads to inefficient use of data as the approach of \cite{sambharya_data-driven_2024} does not allow taking $n_a =  n^*_N(\rho)$, as this would depend on the $N$ samples. Our approach overcomes both limitations by simultaneously learning and certifying \( n^*_N(\rho) \) through scenario optimization, yielding more informative evaluations and consequently higher confidence.




\section{Data-driven performance guarantees} \label{sec:performance}

Our convergence analysis has focused on determining the minimum number of iterations required across all problem instances \eqref{parametric_optimization_problem} to achieve convergence within a specified tolerance. We now shift perspective and develop a framework for bounding performance metrics.

When solving problem \eqref{parametric_optimization_problem}, we assume that the solver converges to a point \( z^*(x) \in \mathcal{Z^*_{\rm local}}(x) \) within a user-defined tolerance. We define the performance metric after \( k \) iterations as:
\(\phi(z^k(x), x) = \text{dist}(z^k(x), z^*(x)),\)
where \(\text{dist}: \mathbb{R}^n \times \mathbb{R}^n \rightarrow \mathbb{R}_{\geq 0}\)  
is any continuous-valued distance measure. If we want to evaluate a performance metric that is different from the convergence metric used by the solver, we assume that \(\mathcal{Z^*_{\rm local}}(x)\) consists of a single point. However, if the performance metric is the convergence metric, this assumption is unnecessary. For discrete-valued distance measures, the results from Section \ref{section:convergence guarantee} can be employed instead.

We characterize the problem of bounding the performance metric \(\phi(z^k(x), x)\) through the robust optimization problem:
\begin{equation}\label{e max}
    e^k = \sup_{x \in \Delta} \;\; \phi(z^k(x), x),
\end{equation}
where \(e^k\) denotes the exact worst-case performance bound. Since the direct computation of \(e^k\) is infeasible due to the intractable parameter space \(\Delta\), we instead solve a scenario-based approximation using \(N\) i.i.d. samples \(\{x_i\}_{i=1}^N\). 
This yields the tractable scenario program:
\begin{equation}\label{e_scenario}
\begin{split}
    & e^k_N =  \min_{e \geq 0} \; e\\
    \text{subject to:} \quad & e \geq \phi(z^k(x_i), x_i), \quad \forall i = 1, \ldots, N,
\end{split}
\end{equation}
which simultaneously learns an empirical performance bound \(e^k_N\), and provides probabilistic guarantees on its validity for the entire support \(\Delta\). We impose the following assumption.
\begin{assumption}[Non-accumulation]\label{ass:non-acc}
Given a fixed \(k\), for any \(e \geq 0\),  
\(\mathbb{P}\{x \in \Delta:~\phi(z^k(x), x) = e\} = 0\).
\end{assumption}

The risk \( V(e^k_N) := \mathbb{P}\{\phi(z^k(x), x) > e^k_N\} \), which specializes Definition \ref{def:risk} to the context, quantifies the probability that future instances will exceed the empirical performance bound \( e^k_N \). The convex problem \eqref{e_scenario} has a unique solution, so Theorem \ref{theorem: robust} applies. Moreover, following \cite{campi2008exact}, problem \eqref{e_scenario} is a fully-supported program under Assumption \ref{ass:non-acc}, implying that the probabilistic guarantee in Theorem \ref{theorem: robust} will be tight. Thus, we have:
\begin{equation}\label{e PAC statement}
\begin{split}
    \mathbb{P}^N \{\{x_i\}^N_{i=1}:~ V(e^k_N) \leq \epsilon\} = 1 - \beta,
\end{split}
\end{equation}
with $\epsilon = 1 - \beta^{\frac{1}{N}}$. 

Bounding the performance metric involves two key trade-offs: (i) the trade-off between computation time and metric bound, revealed by solving \eqref{e_scenario} for different iteration counts \( k \), and (ii) the trade-off between violation probability and metric bound at any fixed iteration count $k$, which can be controlled through relaxation of the robust scenario program \eqref{e_scenario} as follows:
\begin{equation}\label{e_relaxed_scenario}
\begin{split}
    \min_{e \geq 0, \, \xi_{N,i} \geq 0, \, i=1,\ldots,N} &  e + \rho \sum^N_{i = 1} \xi_{N,i} \\
    \text{subject to:} \quad e + \xi_{N,i}  \geq & \, \phi(z^k(x_i), x_i), \quad \forall i = 1, \ldots, N.
\end{split}
\end{equation}

Problem \eqref{e_relaxed_scenario} balances the two competing objectives, that is, minimizing the performance bound \( e \) and the constraint violation. The optimal solution \(e^*_N\)\footnote{We assume that the optimal solution \( e_N^* \) is unique; otherwise, a particular solution can be singled-out using a tie-break rule.} is also denoted \(e^k_N(\rho)\). 

We now draw on the results from \cite{GarattiCampi2022}, whose key equations relating the confidence level $\beta$ to the risk bounds are presented below. For a user-specified confidence level \(\beta \in (0,1)\), consider the following polynomial equation in \( t \) for any \( q = 0, 1, \ldots, N-1 \):
\begin{equation}
\binom{N}{q} t^{N-q} - \frac{\beta}{2N} \sum_{i=q}^{N-1} \binom{i}{q} t^{i-q} - \frac{\beta}{6N} \sum_{i=N+1}^{4N} \binom{i}{q} t^{i-q} = 0,
\label{eq:compute_t_2}
\end{equation}
and for $q = N$ consider the polynomial equation:
\begin{equation}
1 - \frac{\beta}{6N} \sum_{i=N+1}^{4N} \binom{i}{q} t^{i-N} = 0.
\label{eq:compute_e_2}
\end{equation}
For any $q = 0, 1, \ldots, N-1$, equation \eqref{eq:compute_t_2} has exactly two solutions in \([0, +\infty)\), which we denote with \( \underline{t}(q) \) and \( \bar{t}(q) \) (\( \underline{t}(q) \leq \bar{t}(q) \)). Equation \eqref{eq:compute_e_2} has only one solution in \([0, +\infty)\), which we denote by \( \bar{t}(N) \), while we define \( \underline{t}(N) = 0 \). Let \( \underline{\epsilon}(q) := \max\{0, 1 - \bar{t}(q)\} \) and \( \bar{\epsilon}(q) := 1 - \underline{t}(q) \) for \( q = 0, 1, \ldots, N \).

\begin{theorem} Consider the relaxed scenario optimization problem \eqref{e_relaxed_scenario}. Fix $\beta \in (0,1)$. Under Assumption \ref{ass:non-acc}, we have
\begin{equation}\label{e relaxed PAC statement}
\begin{split}
    \mathbb{P}^N \{\{x_i\}^N_{i=1}: \underline{\epsilon}(s^*_N) \leq V(e^k_N(\rho)) \leq \overline{\epsilon}(s^*_N)\} \geq 1 - \beta,
\end{split}
\end{equation}
where $\underline{\epsilon}(q), \bar{\epsilon}(q), q = 0,\ldots,N$ are defined as in \eqref{eq:compute_e_2}, and $s_N^* = 1 + q_N^*$.
\end{theorem}

\textit{Proof.} Problem \eqref{e_relaxed_scenario} is always feasible and has a unique solution. The complexity \(s^*_N = 1 + q_N^*\) follows the same reasoning as in the proof of Theorem 2, where \(q_N^*\) is the number of constraints in \eqref{e_relaxed_scenario} for which the relaxation variables are positive.  \hfill \(\square\)

\textbf{Remark 2.} For specific choices of the performance metric $\phi(z^k(x), x)$ and iteration count $k$, Assumption \ref{ass:non-acc} might not hold. Nevertheless, the generalization guarantee established with Theorem 2 remains valid. 

\section{Numerical examples}
\label{Numerical examples}
In this section, we consider an MPC problem, where the current system state $x(k)$ serves as the problem parameter in \eqref{parametric_optimization_problem}, which is formulated as:
\begin{equation}\label{MPC}
\begin{split}
\min_{u_{i|k}, x_{i|k}}  \; \sum_{i=0}^{T-1} (||x_{i|k}&||^2 + ||u_{i|k}||^2) + V_f(x_{T|k})\\
\text{subject to } & \forall i = 0, 1,...,T-1:\\
& x_{0|k} = x(k), \\
& x_{i+1|k} = Ax_{i|k} + Bu_{i|k},\\
& u_{i|k} \in \mathbb{U},\\
& x_{i|k} \in \mathbb{X},\\
& x_{T|k} \in \mathbb{X}_f.
\end{split}
\end{equation}
%
%
We use $A = \begin{bmatrix}
1.5 & 1\\
0 & 1
\end{bmatrix}, B = \begin{bmatrix}
    1.5\\
    1
\end{bmatrix}$, prediction horizon $T = 10$, with state constraints $||x||_\infty \leq 10$, input constraints $||u||_\infty \leq 5$, terminal state constraints $\mathbb{X}_f = \mathbf{0}$ and terminal cost $V_f(x_{T|k}) = 0$. The optimization variables $u_{i|k}$ and $x_{i|k}$ are collectively denoted as the decision vector $z$ in \eqref{parametric_optimization_problem}. The parameter space $\Delta$ inducing feasible problem instances is unknown, we cannot sample from it directly. Instead, we generate $1000$ i.i.d. samples of $x(k) \in \Delta$ by sampling uniformly from the state constraint set $[-10, 10] \times [-10, 10]$, rejecting those that lead to infeasible MPC problem instances. A sampled instance of the parameter \( x(k) \) is denoted as \( x_i \). This approach is valid as our method applies to any probability distribution $\mathbb{P}$ over $\Delta$. 

\subsection{Data-driven convergence assessment}
We first conducted the assessment following Section \ref{section:convergence guarantee}.B.
For each sample $x_i$, we solve the MPC problem using CVXPY \cite{cvxpy} with the solver OSQP \cite{osqp} and record the required iterations $n(x_i)$. We find an empirical convergence threshold $n^*_N = 8200$. Applying Theorem \ref{theorem: robust} with $\beta = 10^{-4}$ yields $\epsilon = 0.009168$, that is, ensuring $ \mathbb{P}\{x \in \Delta: n(x) > 8200\} \leq 0.9168\%$, with confidence at least $1- 10^{-4}$.

\begin{figure}[ht]
    \centering   \includegraphics[width=0.45\textwidth]{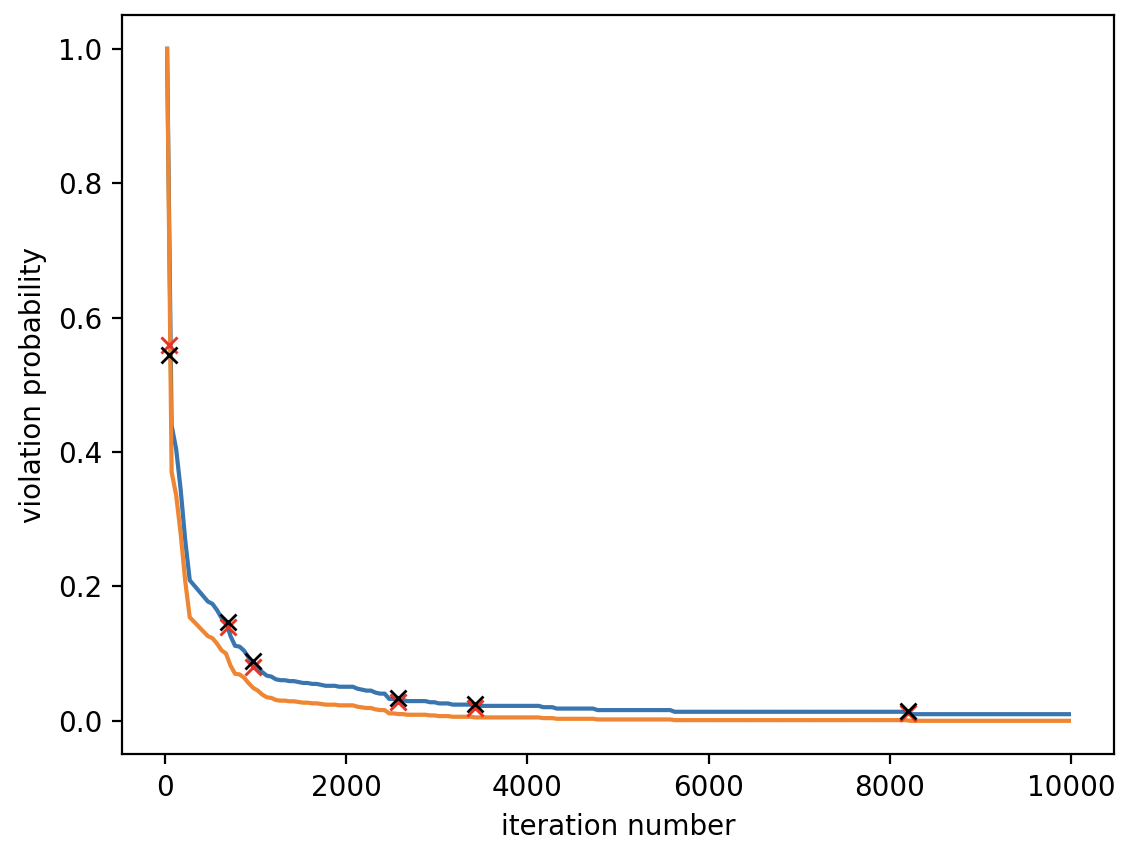}
    \caption{The black crosses indicate the optimal iteration number $n^*_N(\rho)$ and the associated bound $\epsilon(s^*_N)$ computed using Theorem 2, which together provide an upper bound on $V(n^*_N(\rho))$. The blue line shows the upper bound on $V(n_a)$, while the orange line represents the empirical violation probability. The red crosses correspond to the method described in footnote $2$.}
    \vspace{-14pt}
    \label{fig:n trade off}
\end{figure}

We next demonstrate how the relaxed scenario optimization problem \eqref{n_relaxed scenario optimization} can balance the trade-off between the number of iterations and the probability of exceeding these optimal iteration numbers. We follow the process suggested at the end of Section II to select the weight parameter $\rho$.  We first specify the number of constraints we targeted to violate; in this example, we conduct a family of problems, with $\hat{q}_N$ and $\rho$ taking the following values: $\hat{q}_N = 500, 100, 50, 10, 5, 1$, and $\rho = 0.002, 0.01, 0.02, 0.1, 0.2, 1$. We then solve \eqref{n_relaxed scenario optimization} for every $\rho$ to obtain the optimal iteration number $n^*_N(\rho)$ and $q^*_N$, and apply Theorem 2 to compute $\epsilon(s^*_N)$, which is an upper bound for $V(n^*_N(\rho))$. The resulting values of $n^*_N(\rho)$ and $\epsilon(s^*_N)$ are shown with the black crosses in Figure \ref{fig:n trade off}. Note that the value of $\hat{q}_N$ is a good proxy of $q^*_N$ --- in the experiment $q^*_N = 467, 96, 49, 10, 5, 0$, very close to the target values. 

Our framework leverages scenario optimization to bound the risk \( V(n^*_N(\rho)) \), where the iteration count \( n^*_N(\rho) \) is inferred from data. This contrasts with the test set-based method proposed by \cite{sambharya_data-driven_2024}, which bounds \( V(n_a) \) at user-defined $n_a$. Here we use \( n_a = 25, 50, \dots, 9975, 10000 \), represented by the blue curve in Figure \ref{fig:n trade off}. Note that our method naturally focuses on analyzing complex and informative regions, while devoting less attention to flatter, less informative areas. This is clear from the concentration of crosses in regions with steep slopes, and their sparse distribution in flat regions  with small slopes. This dynamical allocation significantly reduces the required evaluations. Consequently, our method achieves a higher overall confidence level, as the confidence level decreases proportionally with the number of evaluations. The orange line indicates the empirically computed violation probability; the fact the theoretical bounds (black crosses) lie very close to this line indicate the tightness of our bound for this case.

\subsection{Data-driven performance assessment}
\begin{figure}[ht]
    \centering   \includegraphics[width=0.45\textwidth]{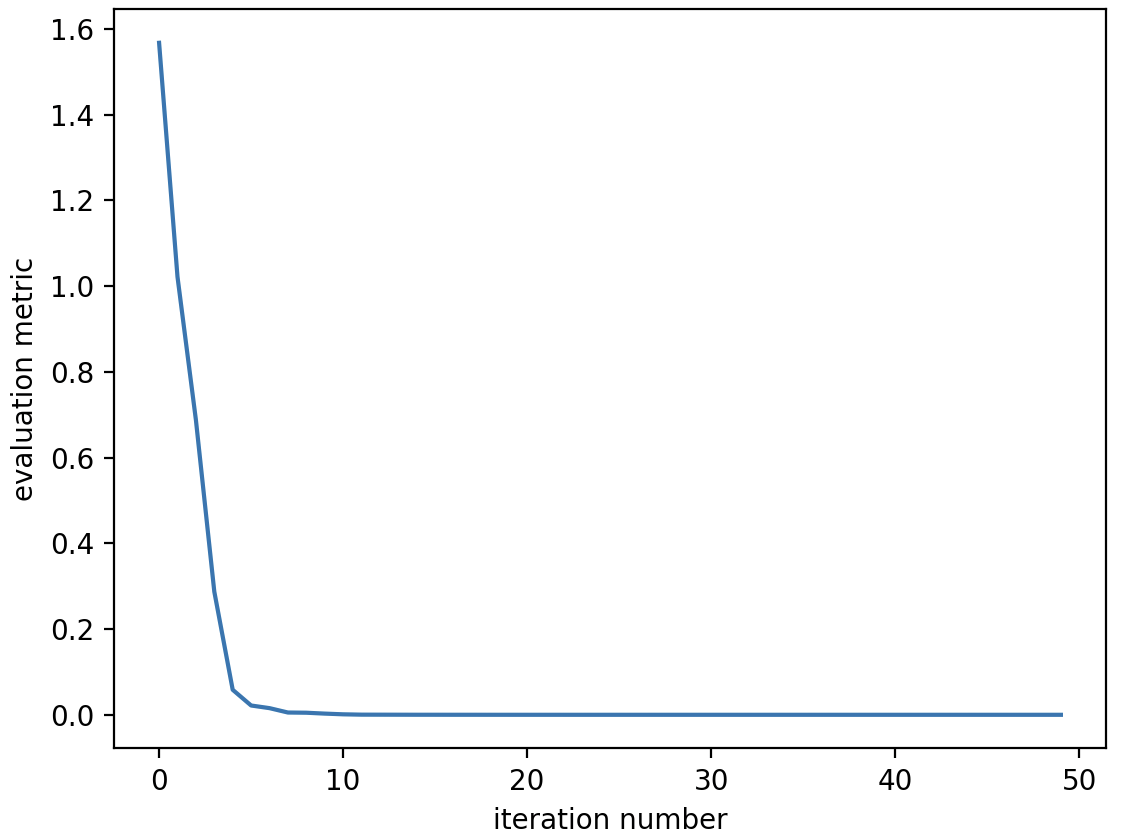}
    \caption{
    Upper bounds on the performance metric computed via \eqref{e_scenario} using $1000$ sampled initial conditions.}
    \vspace{-14pt}
    \label{fig:e_robust}
\end{figure}

We conducted the performance assessment following Section~\ref{sec:performance}. The performance metric considered in this example is $\phi(z^k(x_i), x_i) = \frac{1}{T} \sum_{i=1}^{T} ||x^k_{i|k}-x^*_{i|k}||$. Using the high-precision interior-point solver Clarabel \cite{Clarabel_2024}, we solved the $1000$ sampled instances with permitted iteration counts \(k = 1,...,50\), and record the metrics. Via \eqref{e_scenario}, we establish a relationship between the error bound $e^k_N$ and permitted iteration count $k$, shown in Figure \ref{fig:e_robust}. Each point in this characterization carries the generalization guarantee from \eqref{eq:theorem 1}, ensuring that with confidence at least \( 1-10^{-4} \), the violation probability remains below \( \epsilon = 0.9168\% \).

\begin{figure}[ht]
    \centering   \includegraphics[width=0.45\textwidth]{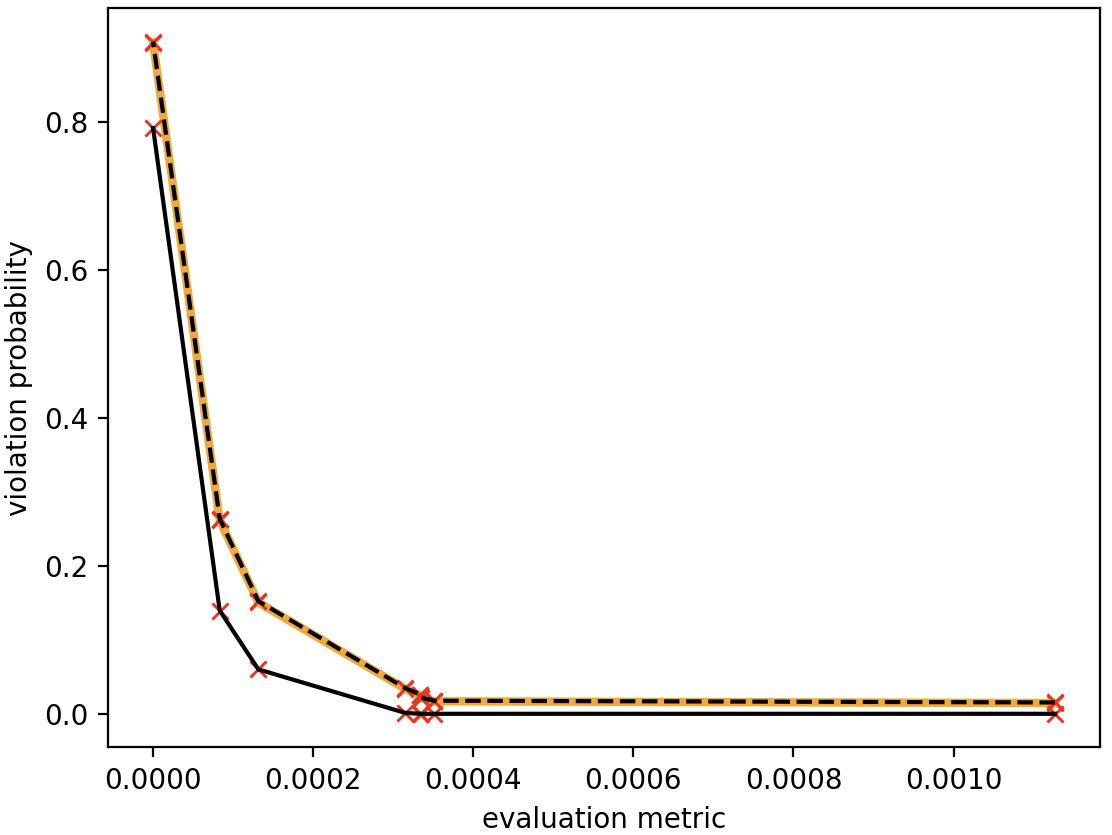}
    \caption{Trade-off between violation probability and metric bound $e^k_N(\rho)$. The black curves represent the bounds $\underline{\epsilon}(s^*_N)$ and $\bar{\epsilon}(s^*_N)$, obtained through Theorem 3, while the orange curve shows bound $\epsilon(s^*_N)$ from Theorem 2.}
    \vspace{-14pt}
    \label{fig:e_relaxed}
\end{figure}

We further implement the formulation with constraint relaxation \eqref{e_relaxed_scenario} using the same sampled instances with the ones of Figure \ref{fig:e_relaxed}, where every point has a confidence of $1-10^{-4}$. For this analysis, we fix the permitted iteration count at $k=11$ and solve \eqref{e_relaxed_scenario} to obtain both the optimized metric bound $e^k_N(\rho)$ and the empirical violation count $q^*_N$,
which is used to compute the bounds of the risk $V(e^k_N(\rho))$. The figure illustrates the trade-off between the violation probability and the metric bound $e^k_N(\rho)$. The black curves represent the upper and lower bounds $\bar{\epsilon}(s^*_N)$ and $\underline{\epsilon}(s^*_N)$, obtained through Theorem 3, which utilizes the non-accumulation assumption, while the orange curve shows bound $\epsilon(s^*_N)$ from Theorem 2 that does not require this assumption. Notably, the orange curve merely coincides with the upper black bound, demonstrating that the non-accumulation assumption provides minimal improvement for the upper bound. However, when valid, this assumption does yield a useful lower bound. 
%
%

\section{Conclusions}
We have proposed a data-driven framework to provide performance guarantees for parametric optimization problems. We consider two types of performance guarantees. The first is how many iterations are required to achieve a predefined tolerance, while the second is an upper bound on a performance metric after a fixed number of iterations. We formulate these two problems as robust scenario optimization problems so that we can obtain empirical estimates of the two performance metrics, alongside probabilistic guarantees for how well those empirical values generalize to unseen problem instances sampled from the same distribution. We then relax the robust scenario optimization problems, and apply the latest advances in the scenario approach to evaluate the trade-off between the solution accuracy and computational cost. Our approach allows practitioners to set iteration budgets that balance computational efficiency and convergence performance in online optimization.

\section*{Acknowledgments}
The authors would like to thank Dr. Guillaume Berger for bringing \cite{New_SVM} to their attention and for suggesting how the a priori method presented in Theorem 1 could be used to provide probabilistic guarantees for problem~\eqref{n_relaxed scenario optimization}.


\bibliographystyle{IEEEtran}
\bibliography{mybib.bib}

\begin{thebibliography}{10}
\providecommand{\url}[1]{#1}
\csname url@samestyle\endcsname
\providecommand{\newblock}{\relax}
\providecommand{\bibinfo}[2]{#2}
\providecommand{\BIBentrySTDinterwordspacing}{\spaceskip=0pt\relax}
\providecommand{\BIBentryALTinterwordstretchfactor}{4}
\providecommand{\BIBentryALTinterwordspacing}{\spaceskip=\fontdimen2\font plus
\BIBentryALTinterwordstretchfactor\fontdimen3\font minus \fontdimen4\font\relax}
\providecommand{\BIBforeignlanguage}[2]{{%
\expandafter\ifx\csname l@#1\endcsname\relax
\typeout{** WARNING: IEEEtran.bst: No hyphenation pattern has been}%
\typeout{** loaded for the language `#1'. Using the pattern for}%
\typeout{** the default language instead.}%
\else
\language=\csname l@#1\endcsname
\fi
#2}}
\providecommand{\BIBdecl}{\relax}
\BIBdecl

\bibitem{MPC_book}
J.~B. Rawlings, D.~Q. Mayne, M.~Diehl \emph{et~al.}, \emph{Model predictive control: theory, computation, and design}.\hskip 1em plus 0.5em minus 0.4em\relax Nob Hill Publishing Madison, WI, 2017, vol.~2.

\bibitem{maciejowski2007predictive}
J.~M. Maciejowski and M.~Huzmezan, ``Predictive control,'' in \emph{Robust Flight Control: A Design Challenge}.\hskip 1em plus 0.5em minus 0.4em\relax Springer, 2007, pp. 125--134.

\bibitem{Jingyi_MPC}
J.~Huang, H.~Wang, K.~Margellos, and P.~Goulart, ``Predictive control barrier functions: Bridging model predictive control and control barrier functions,'' in \emph{Eur. Control Conf. (ECC)}, 2025, pp. 2623--2629.

\bibitem{Beck2017}
A.~Beck, \emph{First-order Methods in Optimization}.\hskip 1em plus 0.5em minus 0.4em\relax Society for Industrial and Applied Mathematics, 2017.

\bibitem{RyuYin2022}
E.~K. Ryu and W.~Yin, \emph{Large-Scale Convex Optimization: Algorithms \& Analyses via Monotone Operators}.\hskip 1em plus 0.5em minus 0.4em\relax Cambridge University Press, 2022.

\bibitem{DasGupta2024}
S.~D. Gupta, R.~M. Freund, X.~A. Sun \emph{et~al.}, ``Nonlinear conjugate gradient methods: worst-case convergence rates via computer-assisted analyses,'' \emph{Math. Progr.}, 2024.

\bibitem{Taylor2017}
A.~B. Taylor, J.~M. Hendrickx, and F.~Glineur, ``Smooth strongly convex interpolation and exact worst-case performance of first-order methods,'' \emph{Math. Progr.}, vol. 161, pp. 307--345, 2017.

\bibitem{sambharya_data-driven_2024}
R.~Sambharya and B.~Stellato, ``Data-driven performance guarantees for classical and learned optimizers,'' \emph{JMLR}, vol.~26, no. 171, pp. 1--49, 2025.

\bibitem{Langford2001NotBT}
J.~Langford and R.~Caruana, ``{(Not) bounding the true error},'' in \emph{NeurIPS}, 2001.

\bibitem{JMLR:v6:langford05a}
J.~Langford, ``Tutorial on practical prediction theory for classification,'' \emph{JMLR}, vol.~6, no.~10, pp. 273--306, 2005.

\bibitem{modulating_robustness_2013}
S.~Garatti and M.~C. Campi, ``Modulating robustness in control design: Principles and algorithms,'' \emph{IEEE Control Systems Magazine}, vol.~33, no.~2, pp. 36--51, 2013.

\bibitem{GarattiCampi2022}
------, ``Risk and complexity in scenario optimization,'' \emph{Math. Progr.}, vol. 191, pp. 243--279, 2022.

\bibitem{Marco_CDC_2020}
M.~C. Campi and S.~Garatti, ``Scenario optimization with relaxation: a new tool for design and application to machine learning problems,'' in \emph{IEEE Conf. Decis. Control (CDC)}, 2020, pp. 2463--2468.

\bibitem{wait_and_judge}
M.~Campi and S.~Garatti, ``Wait-and-judge scenario optimization,'' \emph{Math. Progr.}, vol. 167, pp. 155--189, 2018.

\bibitem{kostas_2014}
K.~Margellos, P.~Goulart, and J.~Lygeros, ``On the road between robust optimization and the scenario approach for chance constrained optimization problems,'' \emph{IEEE Trans. Autom. Control}, vol.~59, no.~8, pp. 2258--2263, 2014.

\bibitem{Kostas_2015}
K.~Margellos, M.~Prandini, and J.~Lygeros, ``On the connection between compression learning and scenario based single-stage and cascading optimization problems,'' \emph{IEEE Trans. Autom. Control}, vol.~60, no.~10, pp. 2716--2721, 2015.

\bibitem{non-convex_scenario_optimization2024}
S.~Garatti and M.~C. Campi, ``Non-convex scenario optimization,'' \emph{Math. Progr.}, vol. 209, pp. 557--608, 2025.

\bibitem{campi2008exact}
M.~C. Campi and S.~Garatti, ``The exact feasibility of randomized solutions of uncertain convex programs,'' \emph{SIAM J. Optimiz.}, vol.~19, no.~3, pp. 1211--1230, 2008.

\bibitem{New_SVM}
B.~Schölkopf, A.~J. Smola, R.~C. Williamson, and P.~L. Bartlett, ``New support vector algorithms,'' \emph{Neural Computation}, vol.~12, no.~5, pp. 1207--1245, 2000.

\bibitem{cvxpy}
S.~Diamond and S.~Boyd, ``{CVXPY}: {A} {P}ython-embedded modeling language for convex optimization,'' \emph{JMLR}, vol.~17, no.~83, pp. 1--5, 2016.

\bibitem{osqp}
B.~Stellato, G.~Banjac, P.~Goulart, A.~Bemporad, and S.~Boyd, ``{OSQP}: an operator splitting solver for quadratic programs,'' \emph{Math. Progr. Computation}, vol.~12, no.~4, pp. 637--672, 2020.

\bibitem{Clarabel_2024}
P.~J. Goulart and Y.~Chen, ``Clarabel: An interior-point solver for conic programs with quadratic objectives,'' 2024.

\end{thebibliography}

\end{document}